\documentclass[fleqn,10pt]{wlscirep}

\usepackage[utf8]{inputenc}
\usepackage[T1]{fontenc}
\usepackage{enumitem}
\usepackage{amssymb}
\usepackage{setspace}
\usepackage{bm}
\usepackage{amsthm,mathrsfs} 
\usepackage{enumitem}
\usepackage{graphicx,graphics,epsfig}
\usepackage{subcaption}
\usepackage[colorinlistoftodos]{todonotes}
\usepackage{caption} 
\usepackage{url}
\usepackage{verbatim}
\usepackage{booktabs}
\usepackage{tabularx}
\usepackage{multirow}
\usepackage{float}
\usepackage[most]{tcolorbox}

\usepackage{deep-macro1}

\theoremstyle{plain}
\newtheorem{thm}{Theorem}

\newcommand{\bthm}{\begin{thm}}
\newcommand{\ethm}{\end{thm}}
\newcommand{\GF}{GraField}

\newcommand{\bpf}{\begin{proof}}
\newcommand{\epf}{\end{proof}}

\usepackage{arydshln}
\theoremstyle{definition}
\newtheorem{defn}{Definition}
\newtheorem{rem}{Remark}
\newtheorem{motivation}{Motivation}

\graphicspath{ {Figures/} }

\title{Spectral Graph Analysis: A Unified Explanation and Modern Perspectives}

\author[1,*]{Subhadeep Mukhopadhyay}
\author{Kaijun Wang}
\affil{Temple University, Department of Statistical Science, Philadelphia, Pennsylvania, 19122, U.S.A.}
\affil[*]{Correspondence and requests for materials should be addressed to S.M. (email: deep@temple.edu)}

\begin{abstract}
Complex networks or graphs are ubiquitous in sciences and engineering: biological networks, brain networks,  transportation networks, social networks, and the World Wide Web, to name a few. Spectral graph theory provides a set of useful techniques and models for understanding `patterns of interconnectedness' in a graph. Our prime focus in this paper is on the following question: \textit{Is there a unified explanation and description of the fundamental spectral graph methods?} 
There are at least two reasons to be interested in this question. Firstly, to gain a much deeper and refined understanding of the basic foundational principles, and secondly, to derive rich consequences with practical significance for algorithm design. However, despite half a century of research, this question remains one of the most formidable open issues, if not the core problem in modern network science. The achievement of this paper is to take a step towards answering this question by discovering a simple, yet universal statistical logic of spectral graph analysis. The prescribed viewpoint appears to be good enough to accommodate almost all existing spectral graph techniques as a consequence of just one \textit{single} formalism and algorithm.
\vspace{-1em}
\end{abstract}
\begin{document}
\flushbottom
\maketitle
\thispagestyle{empty}

%\noindent Please note: Abbreviations should be introduced at the first mention in the main text – no abbreviations lists. Suggested structure of main text (not enforced) is provided below.

\section{Introduction}\label{sec:01_Intro}
How can we efficiently represent a graph to characterize its intrinsic structure, is the ``Holy-Grail'' question that every graph-analyst wrestles with. Spectral graph theory provides an elegant framework and a formal mathematical language to address this key question, which has produced impressive results across a range of application domains including computer vision, neuroimaging, web search, smart grid, social network, recommender systems, high-dimensional modeling, and many others. Over the past few decades, researchers have made dramatic strides toward developing a suite of spectral graph analysis techniques (e.g., Laplacian, Modularity, Diffusion map, regularized Laplacian, Google PageRank model) with increasing sophistication and specialization. However, no single algorithm can be regarded as a panacea for dealing with the evolving complexities of modern graphs. Therefore, the most important and pressing question for the field today appears to be whether we can develop a unifying language to establish ``bridges'' between a body of heavily-used spectral graph analysis techniques, and thus provide logically connected means for reaching different ends. Undoubtedly, any such formulation would be of great theoretical significance and practical value, that will ultimately provide the applied data scientists clear guidance and a systematic strategy for selecting the proper spectral tools to arrive at a confident conclusion.

To that end, this work attempts to unify the theories of spectral graph analysis, embracing the existing paradigm by purely statistical means. Another benefit of this new viewpoint is that it provides surprising insights into the design of fast and scalable algorithms for large graphs, which are otherwise hard to guess using previous understanding.

\subsection{Spectral Graph Analysis: A Practitioner’s Guide} \label{sec:1.1_sga}
The way spectral graph analysis is currently taught and practiced can be summarized as follows (also known as \emph{spectral heuristics}\cite{spielman2007,von2007}):
\begin{enumerate}[topsep=1.4pt,itemsep=1pt]
  \item Let $\cG=(V,E)$ denotes a  (possibly weighted) undirected graph with a finite set of vertices $|V|=n$, and a set of edges $E$. Represent the graph using weighted adjacency matrix $A$ where $A(x,y;\cG)=w_{xy}$ if the nodes $x$ and $y$ are connected by an edge and $0$ otherwise; weights are non-negative and symmetric. The degrees of the vertices are defined as ${\bm d}=A{\mathbf 1}_n$ and $D = {\rm diag}(d_1,\ldots,d_n) \in \cR^{n \times n}$.
  \item Define ``suitable'' spectral graph matrix (also known as graph ``shift'' operator). Most popular and successful ones are listed below:
  \begin{itemize}[noitemsep,topsep=1pt]
    \item $\cL=D^{-1/2}AD^{-1/2}$  ~~~ Chung (1997)\cite{chung1997}.
    \item $\mathcal{B}=A \,-\,N^{-1} {\bm d} {\bm d}^{T}$ ~~~\,Newman (2006)\cite{newman2006}.
    \item $\calT=D^{-1}A$ ~~~\,Coifman and Lafon (2006)\cite{coifman06}.
    \item Type-I Reg. $\cL_\tau= $ $D_\tau^{-1/2}A\,D_\tau^{-1/2}$\, ~~~\,Chaudhuri et al. (2012)\cite{chaudhuri2012}.
    \item Type-II Reg. $\cL_\tau=D_\tau^{-1/2}\,A_\tau\,D_\tau^{-1/2}$ ~~Amini et al. (2013)\cite{amini2013}.
    \item Google's PageRank $\calT_\al=\al D^{-1}A + (1-\al) F$ ~~Brin and Page (1999)\cite{GOOGLE}.
  \end{itemize} \vspace{-.2em}
 where  $\tau, \al>0$ regularization parameters, $A_\tau=A+(\tau/n) {\bf 1}{\bf 1}^{T}$, $D_\tau$ is a diagonal matrix with elements $A_\tau {\mathbf 1}_n$, $F\in \cR^{n \times n}$ with all entries $1/n$, and $N={\rm Vol}(\cG)=\sum_{x,y}A(x,y)$.
  \item Perform spectral decomposition of the matrix selected at step 2, whose eigenvectors form an orthogonal basis of $\cR^n$. Spectral graph theory seeks to understand the interesting structure of a graph by using leading nontrivial eigenvectors and eigenvalues, first recognized by Fiedler (1973)\cite{fiedler73}, which provide the basis for performing harmonic analysis on graphs.
   \item Compute graph Fourier transform (GFT) by expanding functions defined over the vertices of a graph as a linear combination of the selected eigenbasis (from step 3) and carry out learning tasks such as regression, clustering, classification, smoothing, kriging, etc.
\end{enumerate}
\subsection{Motivating Questions: ``Algorithm of Algorithms''} 
\label{sec:1.2_motiv}
The theory and practice of graph data analysis remain very much unsystematic due to the paucity of attention thus far devoted to their foundational issues; the diversity of poorly related graph learning algorithms attests to this fact. Graph signal processing engineers are often left bewildered by the vastness of the existing literature and hugely diverse zoo of methods (developed by the machine learning community, applied harmonic analysts, physicists, and statisticians). As a further negative repercussion, the existing dry mechanical treatment provides the practitioners with absolutely no clue (other than making blind guesses/trial-and-error approach) on how to adapt the existing theory and algorithms for \emph{yet-unseen} complex graph problems.

This has led to a need to develop a broader perspective on this topic, lest we be guilty of not seeing the forest for the trees. The question of \emph{how} different spectral graph techniques can naturally emerge from some underlying basic principles, plays a key role in clarifying the mystery of \emph{why and when} to use them.  The reality is we still lack general tools and techniques to attack this question in a systematic way. Given the very different character of the existing spectral models, it seems that new kinds of abstractions need to be formulated to address this fundamental challenge. What would these abstractions look like?  All this suggests that a new theory is needed, \emph{not just modification of the old}.

One such promising framework is discussed in this paper, which takes inspiration from the celebrated history of nonparametric spectral analysis of time series \cite{blackman1958,tukey1965data,parzen1967role}. The prescribed viewpoint, as described in Section \ref{sec:02_fundament}, draws on an exciting confluence of modern nonparametric methods, mathematical approximation theory, and Fourier harmonic analysis to unveil a more simple conceptual structure of this subject, avoiding undue complexities. Listed below are highlights of the theoretical and practical significance of this work:
\vskip.4em
\begin{itemize}[topsep=1.4pt,itemsep=1pt]
\item It confers the power to \textit{unify} disparate methods (see Fig. \ref{fig:gds}) of graph data analysis in surprising ways; this is discussed in Sections \ref{sec:algo} - \ref{sec:rspec}.
\item it opens up several possibilities for constructing \textit{specially-designed}  more efficient spectral learning algorithms for complex networks; an example is given in Section \ref{sec:kaijun}, which achieves more than a 350x speedup over conventional methods.
\item Finally, in the broader context, it allows to \textit{bridge the gap} between the two modeling cultures: Statistical (based on nonparametric function approximation and smoothing methods) and Algorithmic (based on matrix theory and numerical linear algebra based techniques)--a job far from trivial.
\end{itemize}
\vskip.2em
{\bf Basic Notation.} For an undirected graph $\cG$ we use the following probabilistic notation. Define network probability mass function by $P(x,y;\cG)=A(x,y;\cG)/$ $\sum_{x,y}A(x,y;\cG)$; Vertex probability mass function by $p(x;\cG)=\sum_{y} P(x,y;\cG)$ with the associated quantile function $Q(u;\cG)$. The network cumulative distribution function will be denoted by $F(x;\cG)=$ $\sum_{j \le x} p(j;\cG)$. Finally, define the important Graph Interaction Function by ${\rm GIF}(x,y;\cG)$$=P(x,y;\cG)/p(x;\cG) p(y;\cG),$ $x,y\in V(\cG)$.

\section{Methods: Fundamentals of Statistical Spectral Graph Analysis}
\label{sec:02_fundament}
\subsection{Graph Correlation Density Field}
\label{sec:2.1_grafield}
Wiener's generalized harmonic analysis \cite{wiener1930generalized} formulation on spectral representation theory of time series \emph{starts} by defining the autocorrelation function (ACF) of a signal. In particular, Wiener-Khinchin theorem asserts ACF and the spectral density are Fourier duals of each other.  Analogous to Wiener's \emph{correlation} technique, we describe a new promising starting point, from which we can develop the whole spectral graph theory systematically by bringing nonparametric function estimation and harmonic analysis perspectives. Both statistical and probabilistic motivations will be given.

\begin{defn} (GraField) \label{def:1_grafield}
For given discrete graph $\cG$ of size $n$, the piecewise-constant bivariate kernel function $\cd: [0,1]^2 \rightarrow \cR_+ \cup \{0\}$, called GraField,
is defined almost everywhere through
\beq \label{eq:gcddef}
\cd(u,v;\cG)\,=\,{\rm GIF}\big[Q(u;\cG),Q(v;\cG)\big]\,=\,\dfrac{P\big( Q(u;\cG), Q(v;\cG)\big)}{p\big( Q(u;\cG) \big) p\big(  Q(v;\cG) \big)},~~~0<u,v<1.
\eeq
\end{defn}

\begin{thm} \label{lemma1} GraField defined in Eq. (\ref{eq:gcddef}) is a positive piecewise-constant kernel satisfying
\[\iint_{[0,1]^2} \cd(u,v;\cG) \dd u \dd v\,~ = \sum_{(i,j)\in \{1,\ldots,n\}^2} \iint_{I_{ij}} \cd(u,v;\cG) \dd u \dd v ~= ~1,\]
where
\[I_{ij}(u,v)= \left\{ \begin{array}{ll}
         \,1, ~~& {\rm if}~ (u,v) \in \left(F(i-1;X), F(i;X)\right] \times \left( F(j-1;Y), F(j;Y)\right]\\
         \,0, ~~& \mbox{elsewhere}.\end{array} \right.\]
\end{thm}

The bivariate step-like shape of the \texttt{GraField} kernel is governed by the (piecewise-constant left continuous) quantile functions $Q(u;\cG)$ and $Q(v;\cG)$ of the \emph{discrete} vertex probability measures. As a result, in the continuum limit (let $(\cG_n)_{n\ge0}$ be a sequence of graphs whose number of vertices tends to infinity $n \rightarrow \infty$), the shape of the piecewise-constant discrete $\cd_n$ approaches a ``continuous field'' over unit interval -- a self-adjoint compact operator on square integrable functions (defined on the graph) with respect to vertex probability measure $p(x;\cG)$ (see Section \ref{sec:diffmap} for more details).

\begin{motivation}
We start with a diffusion based probabilistic interpretation of the GraField kernel. The crucial point to note is that the ``slices'' of the $\cd$ \eqref{eq:gcddef} can be expressed as $p(y|x;\cG)/p(y;\cG)$ in the vertex domain, where the conditional distribution $p(y|x;\cG)=P(y,x;\cG)/p(x;\cG)$, can be interpreted  as transition probability from vertex $x$ to vertex $y$ in one time step. This alternative viewpoint suggests a connection with the random walk on the graph with $p(y;\cG)$ being the stationary probability distribution on the graph, as $\lim_{t\rightarrow \infty} p(t,y|x;\cG)=p(y;\cG)$ regardless of the initial starting point $x$ (moreover, for connected graphs the stationary distribution is unique). Here $p(t,y|x;\cG)$ denotes the probability distribution of a random walk landing at location $y$ at time $t$, starting at the vertex $x$. See Lov\'asz (1993) \cite{lovasz1993random} for an excellent survey on the theory of random walks on graphs.

The GraField function $\cd(u,v;\cG)$ measures how the transition probability $p(y|x;\cG)$ is different from the ``baseline'' stationary distribution (long-run stable behavior) $p(y;\cG)$. That \emph{comparison ratio} is the fundamental interaction function for graphs which we denote by ${\rm GIF}(x,y;\cG)$. This probabilistic interpretation along with Theorem \ref{thm:KL} and \ref{thm:lap} will allow us to integrate the \emph{diffusion map} \cite{coifman06} technique into our general framework.
\end{motivation}

\begin{motivation}
GraField compactly represents the affinity or \emph{strength of ties} (or interactions) between every pair of vertices in the graph.
To make this clear, let us consider the following adjacency matrix of a social network representing $4$ employees of an organization
\[ A=\left( {\large\begin{smallmatrix}
0 ~~ &2 ~ ~ &0  ~~ &0\\[.25em]
2 ~ ~&0  ~~ &3  ~~ &3\\[.25em]
0 ~~ &3  ~~ &0  ~~ &3\\[.25em]
0 ~ ~&3  ~~ &3 ~ ~ &0
\end{smallmatrix}}\right), \]
where the weights reflect numbers of communications (say email messages or coappearances in social events etc.). Our interest lies in understanding the strength of association between the employees, denoted as ${\rm Strength}(x,y)$ for $x,y \in V$. Looking at the matrix $A$ (or equivalently based on the histogram network estimator $p(x,y;\cG)=A/N$ with $N=\sum_{x,y} A(x,y)=22$) one might be tempted to conclude that the link between employee $1$ and $2$ is the weakest, as they have communicated only twice, whereas employees $2,3,$ and $4$ constitute strong ties, as they have interacted more frequently. Now, the surprising fact is that (i) ${\rm Strength}(1,2)$ is twice that of  ${\rm Strength}(2,3)$ and  ${\rm Strength}(2,4)$; also (ii) ${\rm Strength}(1,2)$ is $1.5$ times of ${\rm Strength}(3,4)$! To understand this paradox, compute the empirical GraField kernel matrix (Definition 1) with $(x,y)$th entry $N {\cdot} A(x,y;\cG)/d(x)d(y)$
\[ \cd=\left( {\large\begin{smallmatrix}
0  ~&22/8 ~  &0  ~ &0\\[.25em]
22/8  ~&0 ~  &22/16  ~ &22/16\\[.25em]
0 ~ &22/16  ~ &0  ~ &22/12\\[.25em]
0 ~ &22/16 ~ &22/12 ~  &0
\end{smallmatrix}}\right). \]
This toy example is in fact a small portion (with members $1,9,31$ and $33$) of the famous Zachary's karate club data \cite{zachary}, where the first two members were from  Mr. Hi’s group and the remaining two were from John’s group. The purpose of this illustrative example is not to completely dismiss the adjacency or empirical graphon \cite{lovasz2006} based analysis but to caution the practitioners so as not to confuse the terminology ``strength of association'' with ``weights'' of the adjacency matrix -- the two are very different objects. Existing literature uses them interchangeably without paying much attention. The crux of the matter is: association does \emph{not} depend on the raw edge-density, it is a ``comparison edge-density'' that is captured by the GraField.
\end{motivation}

\begin{motivation}
GraField can also be viewed as properly ``renormalized Graphon,'' which is reminiscent of Wassily Hoeffding's ``standardized distributions'' idea \cite{Hoeff40}. Thus, it can be interpreted as a discrete analogue of copula (the Latin word copula means ``a link, tie, bond'') density for random graphs that captures the underlying correlation field. We study the structure of graphs in the spectral domain via this fundamental graph kernel $\cd$ that characterizes the implicit \emph{connectedness} or \emph{tie-strength} between pairs of vertices.
\end{motivation}

Fourier-type spectral expansion results of the density matrix $\cd$ are discussed in the ensuing section, which is at the heart of our approach. We will demonstrate that this correlation density operator-based formalism provides a useful perspective for spectral analysis of graphs that allows both unification and extension.

\subsection{\KL Representation of Graph} 
\label{sec:2.2}
We define the \KL (KL) representation of a graph $\cG$ based on the spectral expansion of its \texttt{GraField} function $\cd(u,v;\cG)$. Schmidt decomposition \cite{schmidt1907} of $\cd$ yields the following spectral representation theorem of the graph.
\begin{thm} \label{thm:KL} The square integrable graph correlation density kernel $\cd:[0,1]^2 \rightarrow \cR_+ \cup \{0\}$ of two-variables admits the following canonical representation
\beq \label{eq;KLcd}
\cd(u,v;\cG_n)\,=\, 1 \,+\, \sum_{k=1}^{n-1} \la_k \phi_k(u) \phi_k(v),
\eeq
where the non-negative $\la_1 \ge \la_2\ge \cdots \la_{n-1}\ge 0$ are singular values and $\{\phi_k\}_{k\ge 1}$ are the orthonormal singular functions $\langle \phi_j,\phi_k\rangle_{\sL^2[0,1]}=\delta_{jk},$ for $j,k =1,\ldots,n-1$, which can be evaluated as the solution of the following integral equation relation
\beq \label{eq:IE}
\int_{[0,1]} \left[\cd(u,v;\cG_n) -1\right] \phi_k(v) \dd v\,=\,\la_k \phi_k(u),~~~~k=1,2,\ldots,n-1.
\eeq
%where ${\rm T}_{\cd}: L^2[0,1] \rightarrow L^2[0,1]$ is a self-adjoint compact operator on square integrable functions (defined on the graph) with respect to vertex probability measure $p(x;\cG)$.
\end{thm}

\begin{rem} By virtue of the properties of KL-expansion \cite{loeve1955book}, the eigenfunction basis $\phi_k$ satisfying \eqref{eq:IE} provides the optimal low-rank representation of a graph in the mean square error sense. In other words, $\{\phi_k\}$ bases capture the graph topology in the \emph{smallest} embedding dimension and thus carries practical significance for graph compression (this will be more evident in Section \ref{sec:kaijun}). Hence, we can call those functions the \emph{optimal} coordinate functions or Fourier representation bases. Accordingly, the fundamental statistical modeling problem hinges on finding approximate solutions to the optimal graph coordinate system $\{\phi_1,\ldots,\phi_{n-1}\}$ satisfying the integral equation \eqref{eq:IE}.
\end{rem}

\subsection{Nonparametric Spectral Approximation}
\label{sec:2.3_specappr}
We view the \emph{spectral graph learning algorithm} as a method of approximating  $(\la_k,\phi_k)_{k\ge1}$ that satisfies the integral equation \eqref{eq:IE}, corresponding to the graph kernel $\cd(u,v;\cG)$. In practice, often the most important features of a graph can be well characterized and approximated by a few top (dominating) singular-pairs. The statistical estimation problem can be summarized as follows:
\[~~~~A_{n \times n} ~~\mapsto~~\cd ~~\mapsto~~ \left\{ \big(\widehat\la_1,\widehat\phi_1\big),\ldots,  \big(\widehat\la_{n-1},\widehat\phi_{n-1}\big) \right\}~\text{that satisfies Eq.}~\eqref{eq:IE}.~~~~~~\]
This allows to recast the conventional matrix calculus-based graph learning approaches as a functional statistical problem.
\vskip.8em
\noindent\textbf{Projection Methods for Eigenvector Approximation}

The Fourier-type nonparametric spectral approximation method starts by choosing an orthogonal basis $\xi_j$ of the Hilbert space $\sL^2(F;\cG)$ satisfying
\[\sum_x \xi_j(x;\cG) \,\xi_k(x;\cG)\, p(x;\cG)~=~0,~~ \text{for}\, j \ne k.\] 
Define $\eta$-functions (quantile domain unit bases), generated from mother $\xi_j$ by
\[\eta_j(u;\cG)~=~\xi_j\big[ Q(u;\cG)\big],~~~0<u<1,~~\]
as piecewise-constant (left-continuous) functions over the irregular grid $\{0,p(1), p(1)+p(2),\ldots, $ $\sum_{j=1}^n p(j)=1\}$ satisfying $\langle \eta_j, \eta_k\rangle_{\sL^2[0,1]}=0$, if $j\ne k$. Approximate the unknown eigenvectors by the projection, $\mathcal{P}_n\phi_k$, on the  \texttt{span}$\{\eta_j,j=1,\ldots,n\}$ defined by
\beq \label{eq:lc} \phi_k(u) \,\approx\, \mathcal{P}_n\phi_k\,=\, \sum_{j=1}^n \te_{jk}\, \eta_j(u),~~~0<u<1 \eeq
where $\te_{jk}$ are the unknown coefficients to be estimated.

\begin{defn}[GraField density matrix]
\label{def:2_grafmat}
We introduce a generalized concept of matrices associated with graphs called the \texttt{G-matrix}. Define discrete graph transform with respect to an orthonormal system $\eta$ as
\beq \label{eq:ogt} \mathcal{M}[j,k; \eta, \cG] ~= ~\Big\langle \eta_j, \int_0^1(\cd-1)\eta_k \Big\rangle_{\sL^2[0,1]}~~~ \text{for}~ j,k=1,\ldots,n.\eeq
Equivalently, we can define the discrete graph transform to be the coefficient matrix of the orthogonal series expansion of the \texttt{GraField} kernel $\cd(u,v;\cG)$ with respect to the product bases $\{\eta_j \otimes \eta_k\}_{1 \le j,k \le n}$. As a practical significance, this generalization provides a systematic recipe for converting the graph problem into a ``suitable'' matrix problem:
\[~~~\cG(V,E)~~~ \longrightarrow~~~ A_{n\times n}~~~~ \longrightarrow ~ ~~~\cd(u,v;\cG) ~~~~ \xrightarrow[\,\,{\rm Eq.}\, (\ref{eq:ogt})]{\{\eta_1,\ldots,\eta_n\}} ~ ~~~ \mathcal{M}(\eta,\cG)\,\in\, \cR^{n\times n}.~~~~~~\]
\end{defn}

\begin{thm}  \label{thm:GM}
The \texttt{G-matrix} \eqref{eq:ogt} can also be interpreted as a ``covariance'' operator for a discrete graph by recognizing the following equivalent representation 
{\small \[ \mathcal{M}[j,k;\eta \mapsto \xi,\cG]=\iint_{[0,1]^2} \big(\cd(u,v;\cG) 
-1 \big)\eta_j(u;\cG) \eta_k (v;\cG)\dd u\dd v=\Cov\left[\xi_j(X;\cG),\xi_k(X;\cG)\right], ~~1 \le j,k \le n.\]}
\end{thm} \vspace{-.5em}
This dual-representation can be proved using the basic quantile mechanics fact that $Q(F(X)) = X$ holds with probability $1$ \cite{parzen79}. Next, we present a general approximation scheme that provides an effective method of discrete graph analysis in the frequency domain. 
\vskip1em
\begin{thm}[Nonparametric spectral approximation] \label{thm:GEN}
The Fourier coefficients $\{\te_{jk}\}$  of the projection estimators \eqref{eq:lc} of the GraField eigenfunctions (eigenvalues and eigenvectors), satisfying the integral equation \eqref{eq:IE},  can be obtained by solving the following generalized matrix eigenvalue problem:
\beq \label{eq:master} \hskip2.1in\mathcal{M} \Te\,=\,S \Te \Delta, \eeq
where $\mathcal{M}_{jl}=\big\langle \eta_j, \int_0^1(\cd-1)\eta_l \big\rangle_{\sL^2[0,1]}$, $\Te_{jl}=\te_{jl}, \Delta_{jl}=\delta_{jl}\la_l$, and $S_{jl}=\big\langle \eta_j, \eta_l \big\rangle_{\sL^2[0,1]}$.
\end{thm}

\begin{proof}
To prove define the residual of the governing equation \eqref{eq:IE} by expanding $\phi_k$ as series expansion \eqref{eq:lc},
\beq \label{eq:er}
R_k(u) ~\equiv~ \sum_j \te_{jk} \Big[ \int_0^1 \big(  \cd(u,v;\cG) -1 \big) \eta_j(v) \dd v \,-\, \la_k  \eta_j(u)   \Big]\, =\, 0.
\eeq
Now for complete and orthonormal $\{\eta_l\}$ requiring the error $R_k(u)$ to be zero is equivalent to the statement that $R(u)$ is orthogonal to each of the basis functions
\beq \big\langle R_k(u),\eta_l(u)\big\rangle_{\sL^2[0,1]}~=~0,~~~~~~k=1,\ldots,n.  \eeq
This leads to the following set of equations:
\beq \label{eq:m1}
\sum_j \te_{jk} \Big[ \iint\limits_{[0,1]^2} \big(  \cd(u,v;\cG) -1 \big) \eta_j(v) \eta_l(u) \dd v \dd u\Big]~-~\la_l \sum_j \te_{jl} \big[\int_0^1 \eta_j(u) \eta_l(u) \dd u \big]~=~0.
\eeq
Complete the proof by writing \eqref{eq:m1} in a compact matrix format: $\mathcal{M} \Te\,=\,S \Te \Delta $.
\end{proof}

\begin{figure}[t]
\centering
\vspace{-.24em}
\includegraphics[width=\textwidth,trim=1.4cm 15cm 1.4cm 4.5cm]{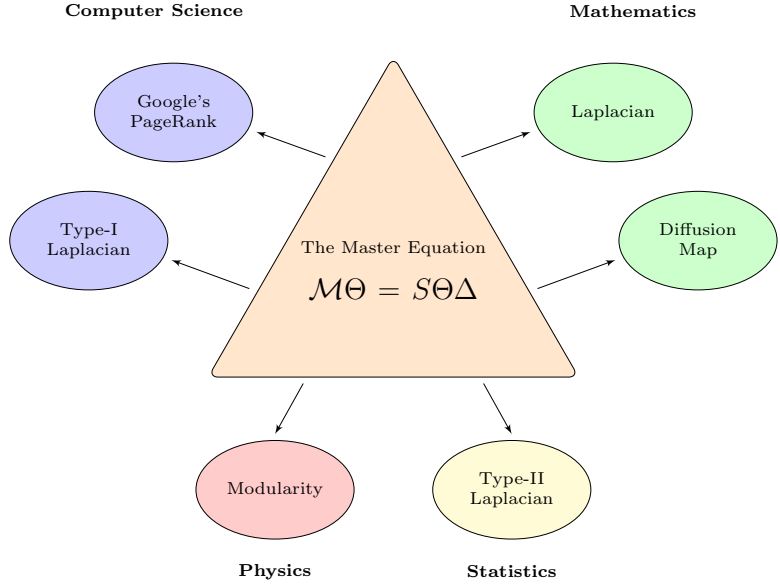}
\caption{The uniting power of our theory that integrates wide-range of spectral graph methods invented by various researchers from different disciplines.} \label{fig:gds}
\vspace{-.4em}
\end{figure}

Few notes on the significances of Theorem \ref{thm:GEN}:
\begin{itemize}
\item The fundamental idea behind the Rietz-Galerkin \cite{galerkin1915} style approximation scheme for solving variational problems in Hilbert space played a pivotal inspiring role to formalize the statistical basis of the proposed computational approach.
\vskip.4em
\item Theorem \ref{thm:GEN} plays a key role in our statistical reformulation. In particular, we will show how the fundamental equation \eqref{eq:master} provides the desired unity among different spectral graph techniques by systematically constructing a ``suitable'' class of coordinate functions.
\vskip.4em
\item Our nonparametric spectral approximation theory based on eigenanalysis of \texttt{G-Matrix}, remains \emph{unchanged} for any choice $\eta$-function, which paves the way for the generalized harmonic analysis of graphs.
\vskip.4em
\item In statistical terms, the search for the new spectral methods amounts to different (wise) choices of $\eta$-functions, and nonparametric estimation methods of GraField $\cd_n$. This systematic view on extending the capability of traditional graph analysis methods (cf. Section \ref{sec:kaijun} for an example) is one of the amusing advantages of our viewpoint, which allows us to discover many ``mysterious similarities'' among classical and modern spectral graph analysis algorithms (Fig \ref{fig:gds} gives the schematic description). This is carried out in Sections \ref{sec:ESM}-\ref{sec:rspec}. 
\end{itemize}
\vskip.4em

\section{Unified Computing Algorithm} \label{sec:algo}
Here we provide the universal construction principle whose significance stems from two considerations: (i) \textit{All} known spectral graph techniques (shown in Figure \ref{fig:gds}) are just different manifestations of the following single general algorithm, as shown in Sections \ref{sec:ESM}-\ref{sec:rspec}. (ii) What is even more surprising is the fact that this \textit{same} general recipe can also be used to design faster, more accurate \textit{custom-made} graph-mining algorithms, as described in Section \ref{sec:kaijun}. 
\vskip1.25em
\begin{center}
\emph{Algorithm 1. A Unified Algorithm of Spectral Graph Analysis} \label{App:Algo1}
\end{center}
\medskip\hrule height .65pt
\vskip1em
% \texttt{Step 1.} For given discrete graph $\cG$ of size $n$, construct GraField kernel function $\cd: [0,1]^2 \rightarrow \cR_+ \cup \{0\}$ defined a.e by
% \beq 
% \cd(u,v;\cG)\,=\,{\rm GIF}\big[Q(u;X),Q(v;Y);\cG\big]\,=\,\dfrac{P\big( Q(u;X), Q(v;Y);\cG \big)}{p\big( Q(u;X) \big) p\big(  Q(v;Y) \big)},
% \eeq
% where $u=F(x;\cG),v=F(y;\cG)$ for $x,y \in \{1,2,\ldots,n\}$.
% \vskip.24em
\texttt{Step 1.}  Let $\{\xi_j(x;\cG)\}_{1\le j \le n}$ to be a discrete orthonormal basis of $\sL^2(F)$ Hilbert space (The vertex probability measure) 
satisfying the following conditions:
\[\sum\nolimits_x \xi_j(x;\cG) \,\xi_k(x;\cG)\, p(x;\cG) \,= \,\delta_{jk}.\]

\noindent\texttt{Step 2.} $\xi \mapsto \eta$ via quantile transform: Construct $\eta_j(u;\cG):= \xi_j(Q(u;\cG))$ on the unit interval $0\le u \le1$.

\vskip.34em

\noindent\texttt{Step 3.} Transform coding of graphs. Construct generalized spectral graph matrix or \texttt{G-matrix} $\mathcal{M}(\xi;\cG)\in \cR^{n\times n}$ by:
\beq \label{eq:ogt1} \mathcal{M}[j,k;\xi,\cG]\,=\Big\langle \eta_j, \int_0^1(\cd-1)\eta_k \Big\rangle_{L^2[0,1]}=\,\sum_{\ell,m} \xi_j(\ell;\cG) \xi_k(m;\cG) P(\ell,m;\cG).\eeq
$\mathcal{M}(\xi;\cG)$ can be viewed as a \emph{transform coefficient matrix} of the orthogonal series expansion of $\cd(u,v;\cG)$ with respect to the product bases $\{\eta_j \otimes  \eta_k\}_{1 \le j,k \le n}$.
\vskip.34em
\noindent\texttt{Step 4.} Perform the singular value decomposition (SVD) of $\mathcal{M}(\xi;\cG)=U\Lambda U^{T}$ $= \sum_k u_k \la_k u_k^{T}$, where $u_{ij}$ are the elements of the singular vector of moment matrix $U=(u_1,\ldots,u_n)$, and $\Lambda={\rm diag}(\la_1,\ldots,\la_n)$, $\la_1 \ge$ $ \cdots \la_n \ge 0$.
\vskip.34em
\noindent\texttt{Step 5.} Obtain approximate \KL basis (or simply the graph representation basis) of the graph $\cG$ by taking the following linear combination:
\[\widetilde{\phi}_{k}\,=\,\sum_{j=1}^n u_{jk}\xi_{j},\,~ {\rm for}\,\, k=1,\ldots,n-1,\]
which can be directly used for subsequent signal processing on graphs.
\vskip.4em
\medskip\hrule height .65pt
\vskip1.5em
\section{Empirical Spectral Graph Analysis} \label{sec:ESM}
Three popular traditional spectral graph analysis models will be synthesized in this section as a special case of the master equation (Eq. \ref{eq:master} of Theorem \ref{thm:GEN}) and the previous algorithm.
\subsection{Laplacian Spectral Analysis}

%%%%%

Laplacian is probably the most heavily used spectral graph technique in practice \cite{chung1997, tenenbaum2000science,shi2000,belkin2001laplacian,kelner2011}. Here we will demonstrate for the first time how the Laplacian of a graph \emph{naturally originates} by purely statistical reasoning, totally free from the classical combinatorial based logic.
\vskip.4em

{\bf Degree-Adaptive Block-pulse Basis Functions}. One of the fundamental, yet universally valid (for any graph) choice for $\{\eta_j\}_{1\le j \le n}$ is the indicator top hat functions (also known as  block-pulse basis functions, or in short BPFs). However, instead of defining the BPFs on a uniform grid (which is the usual practice) here (following Sec \ref{sec:2.3_specappr}) we define them on the non-uniform mesh $0=u_0 < u_1\cdots< u_n=1$ over [0,1], where $u_j=\sum_{x \le j}p(x;X)$ with local support
\beq \label{eq:basis}
\eta_j(u) = \left\{ \begin{array}{ll}
         p^{-1/2}(j) ~~& \mbox{for $u_{j-1} < u \le u_{j}$};\\
         0 ~~& \mbox{elsewhere}.\end{array} \right.
\eeq
They are disjoint, orthogonal, and a complete set of functions satisfying
\[\int_0^1 \eta_j(u) \dd u =\sqrt{p(j)},~~\,\int_0^1\eta_j^2(u) \dd u =1,~~\text{and}~~\int_0^1\eta_j(u)\eta_k(u) \dd u = \delta_{jk}.\]

\begin{figure*}[t]
\centering
\includegraphics[,width=\textwidth,keepaspectratio]{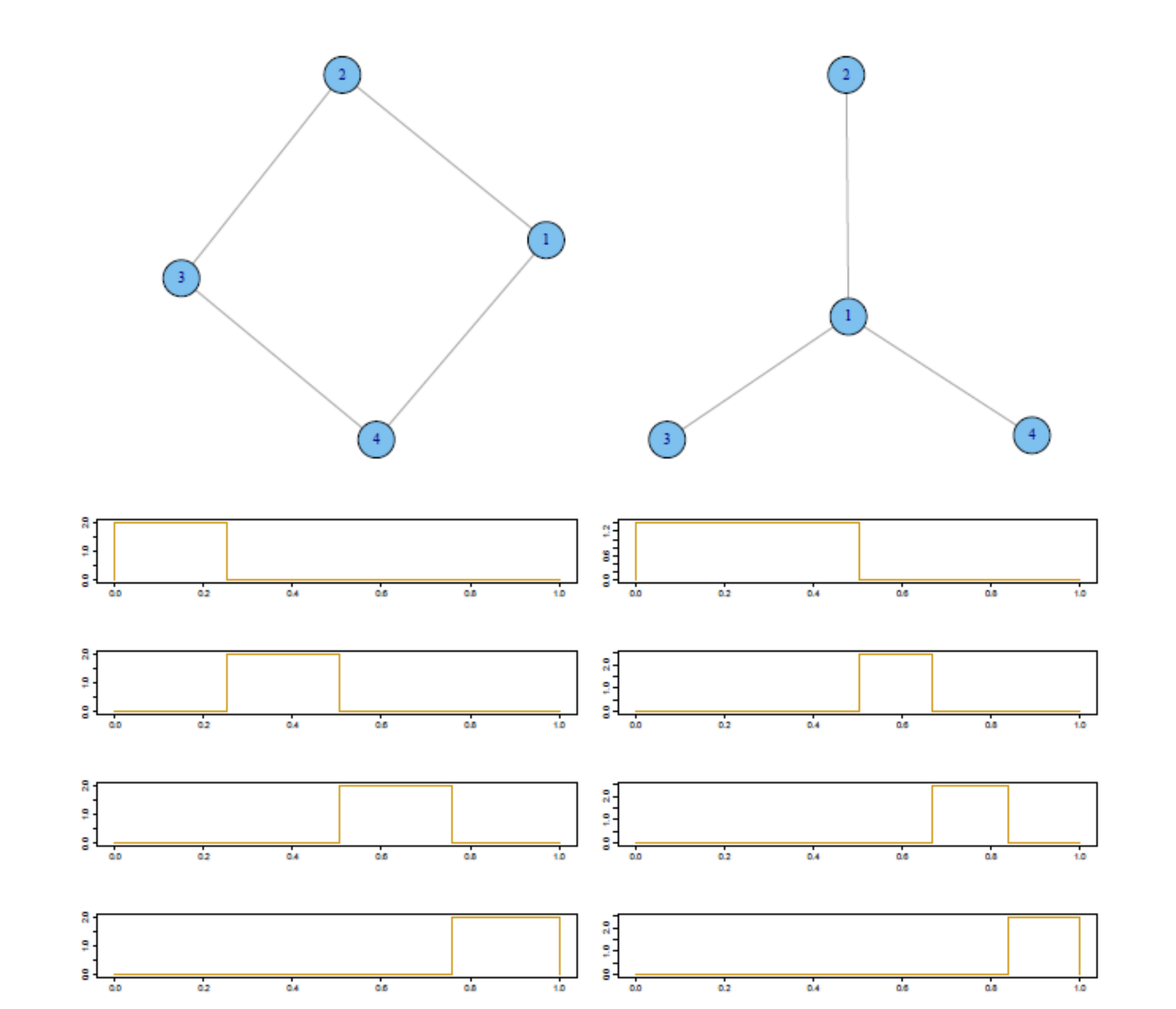}
\caption{Two graphs and their corresponding degree-adaptive block-pulse functions. The amplitudes and the block length of the indicator basis functions \eqref{eq:basis} depends on the degree distribution of that graph.} \label{fig:bpf}
\vspace{-.6em}
\end{figure*}

It is interesting to note that the shape (amplitudes and block lengths) of our specially designed BPFs depend on the specific graph structure via $p(x;\cG)$ as shown in Fig \ref{fig:bpf}. In order to obtain the spectral domain representation of the graph, it is required to estimate the spectra of GraField kernel $\phi_k$, by representing them as block pulse series. The next result describes the required computational scheme for estimating the unknown expansion coefficients $\{\te_{jk}\}$.
\begin{thm} \label{thm:lap}
Let $\phi_1,\ldots,\phi_n$ the canonical Schmidt bases of $\sL^2$ graph kernel $\cd(u,v;\cG)$, satisfying the integral equation \eqref{eq:IE}.
Then the empirical solution of \eqref{eq:IE} for block-pulse orthogonal series approximated \eqref{eq:basis} Fourier coefficients $\{\te_{jk}\}$ can equivalently be written down in closed form as the following matrix eigen-value problem
\beq \label{eq:ev}
\cL^*\Theta\,=\,\Theta \Lambda,
\eeq
where $\cL^*=\cL-uu^{T}$, $\cL$ is the Laplacian matrix, $u=D_p^{1/2}1_n$, and $D_p={\rm diag}(p_1,\ldots,p_n)$.
\end{thm}

\begin{proof}
Recall that the discrete GraField kernel ${\rm GIF}(j,k;\cG)=P(j,k;\cG)/$ $p(j;\cG)p(k;\cG)$ and the tensor product bases $\eta_j(u)\otimes \eta_k(v)$ take the value $p^{-1/2}(j;\cG)\times$ $p^{-1/2}(k;\cG)$ over the rectangle $I_{jk}(u,v)$ for $0<u,v<1$. Substituting this into the master equation \eqref{eq:m1} leads to the following system of linear algebraic equations expressed in the vertex domain:
\beq \label{eq:mm2}
\sum_j \te_{jk} \left[ \dfrac{p(j,l)}{\sqrt{p(j)p(l)}}\,-\,\sqrt{p(j)}\sqrt{p(l)}\,-\,\la_l\delta_{jl} \right] ~=~0.
\eeq
The next step of the proof involves empirical plugin nonparametric estimation as the estimating equation \eqref{eq:mm2} contains unknown network and vertex probability mass functions which need to be estimated from the data. The most basic nonparametric estimates are $\widetilde{P}(j,k;\cG)=A[j,k;\cG]/N$ and $\widetilde{p}(j;\cG)=d(j;\cG)/N$. By plugging these empirical estimators into Eq. \eqref{eq:mm2}, said equation can be rewritten as the following compact matrix form:
\beq \label{eq:m3}
\left[\cL-N^{-1}\sqrt{{\bm d}} \sqrt{{\bm d}}^{T}\right] \tilde{\Theta} = \tilde{\Theta} \tilde{\Lambda} ,
\eeq
where $\cL=D^{-1/2}AD^{-1/2}$ is the graph Laplacian matrix and $\sqrt{{\bm d}}=(\sqrt{ d_1},\ldots,\sqrt{ d_n})^T$. This also leads to an explicit expression for the degree-adaptive BPFs-approximated \KL basis $\widetilde\Phi=D_{\tilde{p}}^{-1/2}\widetilde{\Te}$ by combining \eqref{eq:lc} and \eqref{eq:basis}.
\end{proof}

\vspace{-.5em}
\begin{rem} Theorem \ref{thm:lap} allows us to interpret the graph Laplacian as the empirical \texttt{G-matrix} $\widetilde{\mathcal{M}}(\eta,\cG)$ under degree-adaptive indicator basis choice for the $\eta$ shape function. Our technique provides a completely nonparametric statistical derivation of an algorithmic spectral graph analysis tool. The plot of $\widetilde{\la}_k$ versus $k$ can be considered as a ``raw periodogram'' analogue for graph data analysis.
\end{rem}

\subsection{Diffusion Map} \label{sec:diffmap}
We provide a statistical derivation of Coifman's diffusion map \cite{coifman06} algorithm, which hinges upon the following key result.

\begin{thm} \label{thm:dm}
The \emph{empirical GraField} admits the following vertex-domain spectral \emph{diffusion} decomposition at any finite time $t$
\beq \label{eq:RW} \dfrac{\widetilde p(t,y|x;\cG)}{\widetilde p(y;\cG)}~=~1+\sum_k \tilde \la_k^t \tilde{\phi}_k(x) \tilde{\phi}_k(y),\eeq
where $\widetilde{\phi}_k=  D_{\tilde p}^{-1/2}u_k$, ($u_k$ is the $k$th eigenvector of the Laplacian matrix $\cL$), $\widetilde p(y|x;\cG)=\calT(x,y)$, $\calT=$ $D^{-1}A$ denotes the one-step transition matrix of a random walk on $\cG$ with stationary distribution $\tilde{p}(y;\cG)=d(y;\cG)/N$ and $p(t,y|x;\cG)$ is the kernel of the t-th power of $p(y|x;\cG)$ for $t\ge 0$.
\end{thm}
\begin{proof}
We see that from Theorem \ref{thm:lap} that the Laplacian approximated \GF\, eigen basis is given by $\widetilde{\phi}_k=  D_{\tilde p}^{-1/2}u_k$, where $u_k$ is the $k$th eigenvector of the Laplacian matrix $\cL$. This, together with Theorem \ref{thm:KL}, immediately imply the following vertex-domain spectral decomposition result of the empirical GraField: \beq \label{eq:RW} \dfrac{\widetilde p(y|x;\cG)}{\widetilde p(y;\cG)}~=~1+\sum_k \tilde \la_k \tilde{\phi}_k(x) \tilde{\phi}_k(y).\eeq 
This significance of this identity is that it automatically yields the (appropriately normalized) diffusion kernel whose continuous analogue is known as Laplace-Beltrami operator. As $p(y|x;\cG)$ can be interpreted as the probability of transition in one time step from node $x$ to node $y$, the conclusion is immediate by applying the spectral theory on the $t$-th power of $p(y|x;\cG)$. For an alternative proof of the expansion \eqref{eq:RW} see Appendix section of Coifman
and Lafon (2006) \cite{coifman06} by noting $\widetilde{\phi}_k$ are the right eigenvectors of random walk Laplacian $\calT$.
\end{proof}
\begin{rem}
In light of Theorem \ref{thm:dm}, define diffusion map coordinates at time $t$ as the mapping from $x$ to the vector
\[ x \longmapsto \Big( \la_1^t\phi_1(x), \ldots,  \la_k^t\phi_k(x)\Big),~~~~x\in\{1,2,\ldots,n\},\]
which can be viewed as an approximate `optimal' \KL representation basis and thus can be used for non-linear embedding of graphs. Define diffusion distance, a measure of similarity between two nodes of a graph, as the Euclidean distance in the diffusion map space
\beq \label{eq:DD} D_t^2(x,x')~=~\sum_{j\ge 1} \la_j^{2t}\big\{\phi_j(x)  \,-\,\phi_j(x')\big\}^2.~~~~~~~~~~~\eeq
This procedure is known as \emph{diffusion map}, which has been extremely successful tool for manifold learning\footnote{Manifold learning: Data-driven learning of the ``appropriate'' coordinates to identify the intrinsic nonlinear structure of high-dimensional data. We claim the concept of GraField allows decoupling of the geometrical aspect from the probability distribution on the manifold.}. Our approach provides an \emph{additional} insight and justification for the (Laplace-Beltrami) diffusion coordinates by interpreting it as the strength of connectivity profile for each vertex, thus establishing a close connection with empirical GraField.
\end{rem}

\subsection{Modularity Spectral Analysis}
\begin{thm} \label{thm:Mod}
To approximate the KL graph basis $\phi_k= \sum_j\te_{jk} \eta_j$, choose $\eta_j(u)= \I(u_{j-1} < u \le u_{j})$ to be the characteristic function satisfying
\[\int_0^1 \eta_j(u) \dd u ~=~ \int_0^1 \eta^2_j(u) \dd u ~=~\,p(j;\cG). \]
Then the corresponding empirically estimated spectral graph equation \eqref{eq:m1} can equivalently be reduced to the following  generalized eigenvalue equation:
\beq
\mathcal{B}\Theta\,=\,D \Theta \Lambda,
\eeq
where the matrix $\mathcal{B}$ is given by $A- \frac{1}{{\rm Vol}(\cG)} {\bm d} {\bm d}^T$, and ${\bm d}=A1_n$ is the degree vector.
\end{thm}

\begin{rem}  The matrix $\mathcal{B}$, known as modularity matrix, was introduced by Newman (2006) \cite{newman2006} from an entirely different motivation. Our analysis reveals that the Laplacian and Modularity based spectral graph analyses are equivalent in the sense that they inherently use the same underlying basis expansion (one is a rescaled version of the other) to approximate the optimal graph bases.
\end{rem}
{\bf Summary of Section \ref{sec:ESM}.} Solutions (eigenfunctions) of the GraField estimating equation based on the \texttt{G-matrix} under the proposed specialized nonparametric approximation scheme provides a systematic and unified framework for spectral graph analysis. As an application of this general formulation, we have shown how one can \emph{synthesize the well-known Laplacian, diffusion map, and modularity spectral algorithms} and view them as ``empirical'' spectral graph analysis methods. It is one of those rare occasions where one can witness the convergence of statistical, algorithmic and geometry-motivated computational models.

\section{Smoothed Spectral Graph Analysis} \label{sec:rspec}
Smoothness is a universal requirement for constructing credible nonparametric estimators. Spectral graph analysis is also no exception.
An improved \emph{smooth version} of raw-empirical spectral graph techniques will be discussed, revealing a simple and straightforward statistical explanation of the \emph{origin} of regularized Laplacian techniques. This, as we shall see, also opens up several possibilities for constructing \emph{new} types of spectral regularization schemes, which seems difficult to guess using previous understanding.

\subsection{High-dimensional Undersampled Regime}
Recall from Theorem \ref{thm:GEN} that the generalized matrix eigenvalue equation \eqref{eq:master} depends on the \emph{unknown} network and vertex probability mass functions. This leads us to the question of estimating the unknown distribution $P=(p_1,p_2,\ldots,p_n)$ (support size = size of the graph = $n$) based on $N$ sample, where $N=\sum_{i=1}^n d_i = 2|E|$. Previously (Theorems \ref{thm:lap}-\ref{thm:Mod}) we have used the \emph{unsmoothed} maximum likelihood estimate (MLE) $\tp(x;\cG)$ to construct our \emph{empirical} spectral approximation algorithms, which is the unique minimum variance unbiased estimator. Under standard asymptotic setup, where the dimension of the parameter space $n$ is fixed and the sample size $N$ tends to infinity, the law of large numbers ensures optimality of $\tp$. As a consequence, empirical spectral analysis techniques are expected to work quite well for dense graphs.

\vskip.55em

\emph{Estimation of probabilities from sparse data}. However, the raw $\tp$ is known to be strictly sub-optimal \cite{witten1991} and unreliable in the high-dimensional sparse-regime where $N/n = O(1)$ (i.e., when parameter dimension and the sample size are comparably large). This situation can easily arise for modern day \emph{large sparse graphs} where the ratio $N/n$ is small and there are many nodes with low degree, as is the case of degree sparsity. The naive MLE estimator can become unacceptably noisy (high variability) due to the huge size and sparse nature of the distribution. In order to reduce the fluctuations of ``spiky'' empirical estimates, some form of ``smoothing'' is necessary. The question remains: How to tackle this high-dimensional discrete probability estimation problem, as this directly impacts the quality of our nonparametric spectral approximation.

Our main purpose in the next section is to describe one such promising technique for smoothing raw-empirical probability estimates, which is flexible and in principle can be applied to any sparse data.
\subsection{Spectral Smoothing} \label{sec:addtau} 
We seek a practical solution for circumventing this problem that lends itself to fast computation. The solution, that is both the simplest and remarkably serviceable, is the Laplace/Additive smoothing \cite{laplace1951} and its variants, which excel in sparse regimes \cite{fienberg1973,witten1991}. The MLE and Laplace estimates of the discrete distribution $p(j;\cG)$ are respectively given by
\vskip.55em

$\text{Raw-empirical MLE estimates}:  ~~~~~~\tp(j;\cG)\,=\,\dfrac{d_j}{N}\,;~~~~~~~~~~~~~~~~~~~~~~~~~~~~~~~~~~~$ \\[.4em]
$~~~~~~~~\,\text{Smooth Laplace estimates}:~~~~~~~~~~~~\hp_{\tau}(j;\cG)\,=\,\dfrac{d_j + \tau}{N+n\tau}$~~~~~~~~~$(j=1,\ldots,n)$.
\vskip.4em
Note that the  smoothed distribution $\hp_{\tau}$ can be expressed as a convex combination of the empirical distribution $\tp$ and the discrete uniform distribution $1/n$
\beq \label{eq:lap} \hp_{\tau}(j;\cG)\,=\, \dfrac{N}{N+n\tau}\,\, \tp(j;\cG)\,+\,  \dfrac{n\tau}{N+n\tau} \, \left(\dfrac{1}{n}\right),~~~\eeq
which provides a Stein-type shrinkage estimator of the unknown probability mass function $p$. The shrinkage significantly reduces the variance, at the expense of slightly increasing the bias.
\vskip.4em

{\bf Choice of $\tau$}. The next issue is how to select the ``flattening constant'' $\tau$. The following choices of $\tau$ are most popular in the literature\cite{fienberg1973}:
\[ \tau = \left\{ \begin{array}{ll}
         1 & \mbox{Laplace estimator};\\
        1/2 & \mbox{Krichevsky–Trofimov estimator}; \\
        1/n & \mbox{Perks estimator};\\
       \sqrt{N}/n &\mbox{Minimax estimator (under $L^2$ loss)}.\end{array} \right. \]
\vskip.25em

\begin{rem}
Under increasing-dimension asymptotics, this class of estimator is often difficult to improve without imposing \emph{additional} smoothness constraints on the vertex probabilities. The latter may not be a valid assumption as nodes of a graph offer \emph{no natural order} in general.
\end{rem}

With this understanding, smooth generalizations of empirical spectral graph techniques will be discussed, which have a close connection with recently proposed spectral regularized techniques.

\subsection{Type-I Regularized Graph Laplacian}
Construct $\tau$-regularized smoothed empirical $\eta_{j;\tau}$ basis function by replacing the amplitude $p^{-1/2}(j)$ by $\hat{p}_\tau^{-1/2}(j)$ following \eqref{eq:lap}. Incorporating this regularized trial basis, we have the following modified \texttt{G-matrix} based linear algebraic estimating equation \eqref{eq:master}:
\beq \label{eq:m2}
\sum_j \te_{jk} \left[ \dfrac{\widetilde{p}(j,k)}{\sqrt{\hat p_\tau(j) \hat p_\tau(k)}}\,-\,\sqrt{\hat p_\tau(j)}\sqrt{\hat p_\tau(k)}\,-\,\la_k\delta_{jk} \right] ~=~0.
\eeq
The following theorem follows easily from the preceding calculation.

\begin{thm}
The $\tau$-regularized block-pulse series based spectral approximation scheme is equivalent to representing or embedding discrete graphs in the continuous eigenspace of
\beq
\text{Type-I Regularized Laplacian} ~~=~~D_\tau^{-1/2}\,A\,D_\tau^{-1/2},~~~~~~~~~~~~
\eeq
where $D_\tau$ is a diagonal matrix with $i$-th entry $d_i+\tau$.
\end{thm}
\begin{rem}
It is interesting to note that this exact regularized Laplacian formula was proposed by Chaudhuri et al. (2012) \cite{chaudhuri2012} as well as Qin and Rohe (2013) \cite{qin2013}, albeit from a very different motivation.
\end{rem}

\subsection{Type-II Regularized Graph Laplacian}

\begin{thm} \label{thm:TII}
Estimate the joint probability $p(j,k;\cG)$ by extending the formula given in \eqref{eq:lap} for the two-dimensional case as follows:
\beq \label{eq:2dlap}
\hp_{\tau}(j,k;\cG)\,\,=\,\, \dfrac{N}{N+n\tau}\,\, \tp(j,k;\cG) \,+\,  \dfrac{n\tau}{N+n\tau} \, \left(\dfrac{1}{n^2}\right),~~~~~~
\eeq
which is equivalent to replacing the original adjacency matrix by $A_\tau=A+(\tau/n) {\bf 1}{\bf 1}^{T}$. This modification via smoothing in the estimating equation \eqref{eq:m2} leads to the following spectral graph matrix
\beq \label{eq:rlp2}
\text{Type-II Regularized Laplacian} ~~=~~D_\tau^{-1/2}\,A_\tau\,D_\tau^{-1/2}.~~~~~~~~~~~~
\eeq
\end{thm}
\begin{rem}
Exactly the same form of regularization of Laplacian graph matrix \eqref{eq:rlp2} was proposed by Amini et al. (2013) \cite{amini2013} as a fine-tuned empirical solution.
\end{rem}

\subsection{Google's PageRank Method} \label{sec:google}
Smoothing of network and vertex probability distributions appearing in generalized matrix equation \eqref{eq:master} resulted in Type-I and Type-II regularized Laplacian methods. The third possibility is to directly smooth the conditional or transitional probability matrix to develop a regularized version of random walk Laplacian (which we call Type-III regularized Laplacian) method for large sparse graphs.
\emph{Smoothing conditional probability function}.  Consider a random walk on $\cG$ with transition probability $\calT(i,j;\cG)=\Pr(X_{t+1}=j|X_t=i) \ge 0$. Note that the smoothing \eqref{eq:2dlap} can equivalently be represented as
\beq \label{eq:Google1}
\widehat{\calT}_\tau(i,j;\cG)\,\,=\,\,\dfrac{A(i,j;\cG) + \tau/n}{d_i+\tau}\,\,=\,\, (1-\al_\tau) \widetilde{\calT}(i,j;\cG)\,+\,\al_\tau \left(\dfrac{1}{n}\right),
\eeq
where the degree-adaptive regularization parameter $\al_\tau=\tau(d_i+\tau)^{-1}$ and empirical transition matrix $\widetilde{\calT}=D^{-1}A$.

\begin{thm} \label{thm:Google}
Smooth \textit{each row} of empirical random walk Laplacian $\widetilde{\calT}$ (as it is a row-stochastic or row-Markov matrix) via Stein-like shrinkage to construct an equivalent and more simplified (non-adaptive) estimator:
\beq \label{eq:Google2} \text{Type-III Regularized Laplacian} ~~=~~(1-\al) \widetilde{\calT}\,+\,\al F, ~~~~0 < \al < 1,\eeq
where $F$ is a rank-one matrix given by $\frac{1}{n}{\bf 1}{\bf 1}^T$.
\end{thm}

\begin{rem}
The non-adaptive regularized random-walk Laplacian matrix $\calT_\al$ \eqref{eq:Google2} was introduced by Larry Page and Sergey Brin in 1996 \cite{GOOGLE} and is called the Google's PageRank matrix. What seems natural enough from our nonparametric smoothing perspective, is in fact known to be highly surprising \emph{adjustment}--the ``Teleportation trick''.
\end{rem}
\begin{quote}
`Teleporting is the essential distinguishing feature of the PageRank random walk that had not appeared in the literature before' -- Vigna (2005)\cite{vigna2005}.
\end{quote}

\begin{rem} %Spectral regularization is not a new concept.
The Google's PageRank matrix (which is different from the empirical random walk Laplacian in an important way) is probably the most famous, earliest, and spectacular example of spectral regularization that was originally introduced to counter the problem of \textit{dangling nodes} (nodes with no outgoing edges)  by ``connecting'' the graph\footnote{From a statistical estimation viewpoint, this can be thought of as a way to escape from the ``zero-frequency problem'' for discrete probability estimation.}. As a result the random walk can teleport to a web page uniformly at random (or adaptively based on degree-weights) whenever it hits a dead end. The steady state vector describes the long term visit rate--the PageRank score, computed via eigen-decomposition of $\calT_\al$.
\end{rem}

\subsection{Other Generalizations} \label{sec:ddt}
The beauty of our statistical argument is that it immediately opens up several possibilities to construct \emph{new} types of spectral regularization schemes, which are otherwise hard to guess using previous understanding. Two such promising techniques are discussed here. The first one deals with
Stein smoothing with data-driven shrinkage parameter.
\vskip.4em
We refine the smooth nonparametric spectral graph models (discussed in previous Sections \ref{sec:addtau}-\ref{sec:google}) by wisely choosing the regularization parameter $\tau$: from fixed $\tau$ to data-driven $\widehat\tau$. Quite surprisingly, the explicit formula of the empirical optimal $\widehat \tau^*$ can be derived by minimizing the squared error loss risk function, as given in the following theorem.
\begin{thm} \label{thm:pagerank}
Under the mean square error (MSE) risk function, the Stein-optimal shrinkage parameter and its consistent estimator is given by
\beq \label{final-tau} \widehat \tau^*~=~\dfrac{N^2 - \sum_{i=1}^n d_i^2}{n\sum_{i=1}^n d_i^2 - N^2}.\eeq %N(N+n-1)
\end{thm}
\begin{proof} Our theory of determining optimal shrinkage parameter starts by explicitly writing down the risk function:
\[
{\rm MSE}(\hp_{\tau}, p)\,:=\, R(\tau)\,=\,\sum_{i=1}^n \Ex\big ( \hp_{\tau}(i;\cG) - p(i;\cG) \big)^2,
\]
Replacing by $\hp_{\tau}(i;\cG)$ by \eqref{eq:lap}, after some routine calculations, we get
\beq \label{eq:trisk}
R(\tau)\,=\,\Big( 1- \dfrac{2n\tau}{N+n\tau}\Big) {\rm Var}(\tp(i;\cG))\,+\,\Big(  \dfrac{n\tau}{N+n\tau}\Big)^2 \Ex\big(\tp(i;\cG)-1/n\big)^2,
\eeq
Analytically minimizing the total risk function \eqref{eq:trisk} with respect to $\tau$ by setting $R'(\tau)$ = 0, leads to the following optimal choice:
\beq \label{eq:DDE} \tau^*\,=\,\frac{N}{n} \times \dfrac{\sum_{i=1}^n {\rm Var}(\tp(i;\cG))}{\sum_{i=1}^n \big[ \Ex(\tp(i;\cG) - 1/n)^2 -  {\rm Var}(\tp(i;\cG))\big]}.\eeq
As variance of the MLE estimate ${\rm Var}(\tp(i;\cG))=p(i;\cG)(1-p(i;\cG))/N$, we can further simplify the formula \eqref{eq:DDE} as
\[\tau^*\,=\,\dfrac{1-\sum_{i=1}^n p^2(i;\cG)}{n\sum_{i=1}^n p^2(i;\cG) - 1}.\]
The proof follows by replacing the unknown network probabilities by their empirical counterparts $\tp(i;\cG)=d_i/N$.
\end{proof}

\begin{rem}
The significance of this result lies in the fact that the optimal regularization parameter $\tau$ can be determined with \emph{no} additional computation in a completely data-driven fashion. The formula \eqref{final-tau} provides an analytic expression that guarantees to minimize the MSE avoiding the need of computationally expensive cross-validation like methods. It is shown in Supplementary Appendix D that proposed method performs {\bf 100x} faster and uses {\bf 1000x} less memory for moderately large problems than the current best practice.
\end{rem}

\vspace{-.55em}
While the Laplace or additive smoothing performs well in general, there are situations where they perform poorly \cite{gale1994}.
The Good-Turing estimator \cite{good1953} is often the next best choice which is given by
\beq \label{eq:GT} \widehat{p}_{{\rm GT}}(i;\cG)~=~\dfrac{\varpi_{d_i+1}}{\varpi_{d_i}} \cdot \dfrac{d_i+1}{N},  \eeq
where $\varpi_{d_k}$ denotes the number of nodes with degree $k$. An excellent discussion on this topic can be found in Orlitsky et al. (2003) \cite{orlitsky2003}. One can plug in the estimate \eqref{eq:GT} into the equation \eqref{eq:m2} to generate new spectral graph regularization technique.
\vskip.35em
{\bf Summary of Section \ref{sec:rspec}.} We wish to emphasize that our novelty lies in addressing the open problem of obtaining a \emph{rigorous interpretation and extension} of spectral regularization. To the best of our knowledge this is the \emph{first work} that provides a more formal and intuitive understanding of the \emph{foundation} of spectral regularization. We have shown how the regularization naturally arises as a consequence of high-dimensional discrete data smoothing.

\section{Going Beyond: An Enhanced Method for Accelerated Graph-Learning} \label{sec:kaijun}
Traditional spectral graph analysis techniques do not scale gracefully with the size of the graph, which caused a major practical obstacle to applying them for large-scale problems.  As the size of the datasets is growing rapidly, the key challenge is to find an improved representation of graphs that allows compression. In the following, we demonstrate through examples, how our general theory can be used to find such custom-designed compressors for fast and efficient spectral graph analysis.
\subsection{Motivating Example and Challenges}
We consider the voting records of the US Senate covering the period from the inauguration of George H. W. Bush in January 1989 to the end of Bill Clinton's term in January 2001 (retrieved from www.voteview.com). Table \ref{table:senate} portrays a snapshot of the dataset, where $n=2678$ bills were submitted over the period and binary voting decisions ($1$ for yes, $0$ for no) were recorded from each of the $100$ seats.   
\begin{table}[h]
\caption{Data corresponding to the voting record of the US Senate covering the period from Jan 1989 to Jan 2001. Source: www.voteview.com.} \label{table:senate}
\setlength{\tabcolsep}{11pt}
\renewcommand{\arraystretch}{1.1}
\centering
\begin{tabular}{|c|c|c|c|c|c|c|c|c|}
\hline
\multirow{2}{*}{Time} & \multirow{2}{*}{Bills}&\multicolumn{7}{c|}{Senate Seats}\\
\cline{3-9}
& & ~1~ & ~2~ & ~3~& $\cdots$& 98~ & 99~ & 100\\
\hline
1989-02-28&1& 0 & 1& 1  &$\cdots$& 1& 0  &  1\\
1989-02-28&2&  0 & 1& 0  &$\cdots$& 1 & 0  &  0\\
$\vdots$& &$\vdots$  & &$\vdots$ & & $\vdots$ & & $\vdots$  \\
2000-12-05 &2677& 1 & 1& 1  &$\cdots$& 0& 1  &  1\\
2000-12-07 &2678&  1 & 1& 1  &$\cdots$& 0 & 1  &  1\\
\hline
\end{tabular}
\end{table}
Previously\cite{moody2013portrait,RSSBchangepoint} this dataset was analyzed for studying changes in voting patterns to better understand the political polarization among Republican and Democratic senators in the U.S. Congress. Broadly speaking, this problem falls under the category of 
high-dimensional change-point detection, a critical task in many fields, including econometrics, signal processing, climatology, genomics, neuroscience, speech recognition, and medical monitoring. We observe an ordered sequence of $d$-dimensional data ${\bf Z}_t=(Z_{t1},Z_{t2},\ldots, Z_{td})'$ for $t=1,\ldots, n$ and the interest lies in \textit{tracking} $\{{\bf Z}_t; t=1,\ldots, n\}$ for possible distributional shift:
\[ \hskip1.5in {\bf Z}_t = \left\{ \begin{array}{ll}
         F_1 & \mbox{if $1 \le t \le \tau_1$}\\
        F_2 & \mbox{if $\tau_1 < t \le \tau_2$}\\
        \vdots & ~~~~\vdots\\
         F_k & \mbox{if $\tau_k < t \le \tau_{k+1}=n$}
        ,\end{array} \right. \] 
%which implies the sequence is identically distributed as $F_1$ until a time $\tau_1$, after that the distribution changes abruptly to $F_2$ at $\tau_2$ and so on. Note that here 
where the locations ($\tau_i;i=1,\ldots, k$), the number of change points $k$, and the distributions $F_i,$ $i=1,\ldots,k$ in $\mathbb{R}^d$ are \textit{all} unknown. The goal is simple: to develop a fully nonparametric algorithm for automatic detection of multiple change-points. One way to attack this problem is to transform it into a graph problem. For the Senate data, we construct a network of size $n$ with vertices as the time points (or, equivalently, the bills), and the level of agreement between two voting records was calculated using the Pearson-$\phi^2$ coefficient, a  measure of association for two binary variables. This time-indexed graph, which we call $\mathbb{T}$-graph, encodes the temporal dynamics of the high-dimensional data. In particular, the problem of finding the change-points (splitting the data into different homogeneous time segments) can now be viewed as a spectral clustering exercise on the $\mathbb{T}$-graph. For more discussion on spectral clustering see Supplementary Appendix B. 

As we attempt to partition the network based on $\widetilde\phi_1$, the leading singular vector of the Laplacian matrix as depicted in Fig. \ref{fig:01_sen_kor}, we almost immediately hit a wall. The challenge comes from two principal directions: \textit{statistical and computational}.  The general-purpose Laplacian yields an extremely noisy embedding (the blue dots in Fig. \ref{fig:01_sen_kor}). As a result, the spectral clustering completely fails to segment the time periods into homogeneous blocks to identify the `smooth' transition of the voting patterns around Jan, 1995. The problem is further aggravated by the fact that the classical spectral-partitioning methods are prohibitively expensive, with an $O(n^2)$ cost. In particular, explicitly computing the SVD of the $n \times n$ Laplacian matrix may not be practical or  even feasible for massive-scale problems. Rather than implementing off-the-shelf existing methods in a brute-force manner, one may wonder whether there exists any elegant approach for smoothed and accelerated spectral graph learning. What if we could take advantage of the `internal structure' of the problem? Note that our $\mathbb{T}$-graph is endowed with a special structure, \textit{not} an arbitrary one. The nodes have some kind of natural ordering, which introduces a smoothly varying shape in the Laplacian singular vector, as is clearly visible in Fig. \ref{fig:01_sen_kor}. Once we recognize this, the following question becomes obvious: how to build an efficient and robust computational algorithm that can leverage the special structure of the $\mathbb{T}$-graph? This is where our new statistical point of view comes in handy .
\subsection{Compressive Trial-Basis Design} \label{sec:lpbasis}
The key is to intelligently design the trial-basis functions that take into account the additional structure of the problem. In particular,  $\{\phi_k\}$ must be \textit{compressible} in the new basis, i.e., we need fewer terms (say, $m$, where $m \ll n)$ in the expansion \eqref{eq:lc} to achieve the target accuracy. One such discrete-orthonormal system (called LP-graph basis) is given below, inspired by the recent work of Mukhopadhyay, Parzen, and their colleague\cite{D12e,Deep14LP,Deep18nature}.

To set the stage, we start with a brief overview of the manufacturing process for an arbitrary discrete random variable $X$ with pmf $p(x;X)$ and cdf $F(x;X)$. Construct the basis functions of $X$ by Gram Schmidt orthonormalization of powers of $T_1(x;X)$, given by
\beq T_1(x;X) \,=\,\dfrac{\sqrt{12} [\Fm(x;X) - .5]}{\sqrt{1-\sum_x p^3(x;X)}},\eeq
where $\Fm(x;X)=F(x;X)-.5p(x;X)$ denotes the mid-distribution function. It is easy to verify that these basis functions satisfy the desired orthonormality conditions (see Sec 2.3):
\beq \label{eq:ortho}
\sum_i T_j(x_i;X) p(x_i;X)=0, \sum_i T^2_j(x_i;X) p(x_i;X)=1, \sum_i T_j(x_i;X)  T_k(x_i;X) p(x_i;X) = 0.
\eeq
We define LP orthonormal basis in the \emph{unit interval} (rescaling operation) by evaluating $T_j$ at $Q(u;X)$, quantile function of the random variable $X$
\beq \label{eq:LPbasis} S_j(u;X)\,=\,T_j(Q(u;X);X), ~~0<u<1. \eeq
For undirected weighted graph $\cG$, construct the data-adaptive basis $\{S_j(u;\cG)\}$ by choosing the discrete distribution to be the empirical vertex probability mass function $\tp(x;\cG)=\sum_{y=1}^n A(x,y;\cG)/{\rm Vol}(\cG)$, where ${\rm Vol}(\cG)=\sum_{x,y}A(x,y;\cG)$.  The Eq \eqref{eq:ortho} guarantees that the orthonormal system is automatically \emph{degree-weighted}. One notable aspect of the LP-graph basis is its
nonparametric nature (not fixed like discrete Fourier basis, wavelets, or curvelets, etc.), in the sense that its shape depends on the given graph structure.
\subsection{LP-Frequency Domain Analysis}

To prepare for the main result, here we introduce a Fourier-like transform, called LP-transform coding of graphs, define by
\beq \label{eq:LPT}
\LP[j,k;\cG]\,~~=\,~~\big\langle T_j, T_k \big \rangle_{L^2(\tp)}\,~~=\,~~\big\langle S_j, \int_0^1\cd\,S_k \big\rangle_{L^2[0,1]}~~~~(j,k=1,\ldots,m).
\eeq
This is in fact the \texttt{G-matrix} with respect to an LP-orthonormal system \eqref{eq:ogt1}.
As a practical significance, this provides a systematic recipe for converting the graph problem into a ``suitable'' data-driven reduced dimensional ($m \ll n$)  matrix problem:
\[~~~\cG(V,E)~~~ \longrightarrow~~~ A_{n\times n}~~~~ \longrightarrow ~ ~~~\cd(u,v;\cG) ~~~~ \xrightarrow[\,\,{\rm Eq.}\, (\ref{eq:LPT})]{\{S_1,\ldots,S_m\}} ~ ~~~ \LP\,\in\, \cR^{m\times m}.~~~~~~\]
The purpose of the next result is to describe how to obtain smooth $\widehat \phi_k$ by performing SVD on the smaller-dimensional $\LP_{m \times m}$ matrix instead of operating on the $n \times n$ dimensional Laplacian matrix, thereby significantly accelerating the computation. In other words, how to estimate the unknown coefficients of \eqref{eq:lc} when we approximate $\phi_k(u)$ by $\sum_{j=1}^m \te_{jk} S_j(u)$? The following theorem allows us to pin down the formula.
\begin{thm}\label{thm:LP}
Let $\phi_1,\ldots,\phi_{k_0}$ be the top $k_0$ canonical Schmidt bases of $\sL^2$ graph kernel $\cd(u,v;\cG)$, satisfying the integral equation \eqref{eq:IE}.
Then the empirical solution of \eqref{eq:IE} for LP-orthogonal series approximated \eqref{eq:LPbasis} Fourier coefficients $\{\te_{jk}\}$ can equivalently be written down in closed form as the following matrix eigen-value problem: % in terms of  LP-graph Fourier matrix \eqref{eq:LPT}
\beq \label{eq:ev}
\hskip2in \LP \Theta\,=\,\Theta\Sigma,
\eeq 
which leads to the smooth estimate $\widehat \Phi=S U_{\LP}$, where $\Phi=[\phi_1,\ldots, \phi_{k_0}]$, $S=[S_1,\ldots,S_m]$, and $U_{\LP} \in \cR^{m\times m}$ are the matrix of singular vectors of $\LP$ with singular values $\sigma_1 \ge$ $ \cdots \sigma_m \ge 0$ and $\Sigma={\rm diag}(\sigma_1, \ldots, \sigma_m)$.
\end{thm}
Consequently, we represent the graph in the ``LP-space,'' since important structural information can be compressed into fewer discrete LP-transform coefficients. As an immediate result, we obtain the following reduced-order spectral algorithm. 
% Again, the LP-transformation based nonparametric methods serve as a key vehicle for utilizing the extra structure of the problem.
\begin{center}
\emph{Algorithm 2: LP-Nonparametric Spectral Learning of $\mathbb{T}$-Graph}
\end{center}
\vspace{-.25em}
\medskip\hrule height .65pt
\vskip.5em
\texttt{Step 1.} Inputs: $A$: the weighted adjacency matrix; $m$: LP-smoothing parameter; $k_0$: the desired number of top singular vectors to be approximated.
\vskip.24em
\noindent\texttt{Step 2.} Construct the piecewise-constant orthonormal LP-graph polynomial basis $\{S_1(u;\cG),$ $\ldots, S_m(u;\cG)\}$ using the recipe of Section \ref{sec:lpbasis}:
\vskip.24em
\noindent\texttt{Step 3.} Compute $m \times m$ LP-graph transform matrix by the empirical estimate of \eqref{eq:LPT}:
\[\LP[j,k;\cG]\,=\, \mathop{\sum\sum}_{x,y \in V(\cG)} \tp(x,y;\cG)\, S_j[\wtF(x;\cG)]\, S_k[\wtF(y;\cG)],~~~1 \le j,k\le m.\]
\noindent\texttt{Step 4.} Perform the singular value decomposition (SVD) of $\LP=U_{\LP}\Lambda U_{\LP}^{T}= \sum_k u_k \sigma_k u_k^{T}$, where $u_{ij}$ are the elements of the singular vector of moment matrix $U_{\LP}=(u_1,\ldots,u_m)$, and $\Lambda={\rm diag}(\sigma_1,\ldots,\sigma_m)$, $\sigma_1 \ge$ $ \cdots \sigma_m \ge 0$.
\vskip.24em
\noindent\texttt{Step 5.} Obtain the LP-smoothed graph-Fourier basis by
\[\widehat{\phi}_{k}\,=\,\sum_{j=1}^m u_{jk}S_{j},\,~ {\rm for}\,\, k=1,\ldots,k_0.\]
These ``generalized graph-coordinates'' can now be used directly for subsequent learning.
\vskip.4em
\medskip\hrule height .65pt
\vskip.5em

\begin{figure}[t]
\centering
\includegraphics[width=.65\linewidth]{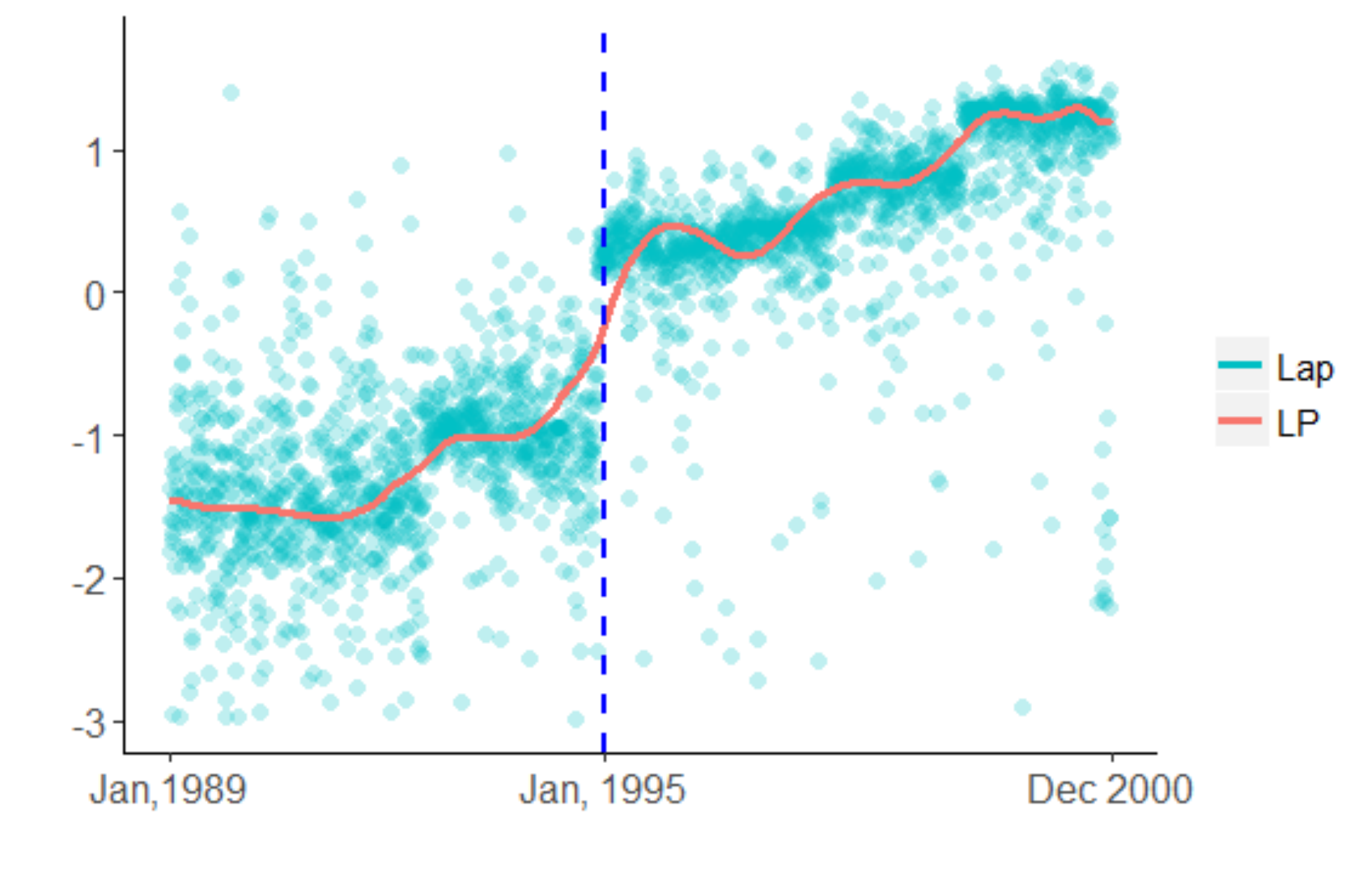}
\vspace{-.4em}
\caption{Spectral plot for Senate data. The `spiky' blue dots denote the Laplacian eigenmap, and the `smooth' red line is the LP-approximated graph-Fourier basis. The dashed blue line indicates the estimated change point (Jan 19, 1995), which is almost indistinguishable from the true change point Jan 17, 1995, when the partisan polarization reached a historically high level.}
\label{fig:01_sen_kor}
\end{figure}

%%%%%%%%%%%%%%%%%%%%%%%%%%%%%%%%%5
\subsection{US Senate Data Analysis}
Our LP-nonparametric spectral approximation algorithm yields the `smooth' red curve in Fig \ref{fig:01_sen_kor}. Contrast this with the noisy Laplacian singular embedding. Applying k-means clustering on the estimated $\widehat \phi_1$ (see Supp Fig. 9 for more details) produces two groups with the change point at January 19th, 1995 (1235th record). This comes close to January 17th, 1995 at the beginning of the tenure of the 104th Congress, when Republicans took control the US House of Representatives for the first time since 1956. According to Moody and Mucha (2013)\cite{moody2013portrait}, this time period witnessed historically high level of partisan polarization not seen in nearly 100 years.

The specially-designed LP-basis functions require only $m$ coefficients (here we have used $m=15$) instead of $n$ to approximate the \KL graph-basis and achieve a remarkable compression ratio of $2678/15=178.53$. This reduces the memory footprint and speeds up the computation. For the Senate data, LP delivers {\bf 350x} speedup compared to the Laplacian-based traditional method! 
%%%%%%%%%%%%%%%%%%%%%%%%%%%%%%%%%%%%%%%%%%%%%%%%%
%%%%%%%%%%%%%%%%%%%%%%%%%%%%%%%%%%%%%%%%%%%%%%%%%

{\bf Summary of Section \ref{sec:kaijun}}. Three unique benefits of our nonparametric spectral algorithm are: compressibility, smoothness, and structure-awareness.  Additional simulation results are reported in Supplementary Appendix E to highlight the robustness and computational efficiency of the proposed scheme.
\section{Concluding Remarks}
Our research takes aim at some foundational questions about complex network modeling that are still unsolved. The prescribed approach brings a fresh statistical perspective to accomplish the miracle of unifying and generalizing the existing paradigm, leading to fundamentally new theoretical and algorithmic developments. It was a great surprise for us to be able to deduce existing techniques (which are often viewed as spectral heuristic-based \emph{empirical mystery facts}) from some underlying basic principles in a self-consistent way. However, from a practical standpoint, we have shown how this new way of thinking about graphs can provide the adequate guidance for designing next-generation computational tools for large-scale problems. It is our hope that a comprehensive understanding gained from this new perspective will provide useful insights for analyzing a wealth of new phenomena that arise in modern network science \cite{marron2007,lazer2009life,kolaczyk2017book,wasserman2016}.

\bibliography{ref-bib}
\vskip.6em
\subsection*{Acknowledgements} The first author benefited greatly from many fruitful discussions and exchanges of ideas he had at the John W. Tukey 100th Birthday Celebration meeting, Princeton University September 18, 2015. This research was also presented at the Parzen Memorial Lecture, 2017 Joint Statistical Meetings--Baltimore, Maryland. We dedicate this research to the beloved memory of Emanuel (Manny) Parzen, 1929--2016.
\subsection*{Author contributions statement} S.M. conceived the research and wrote the paper. K.W. managed the data collection. S.M. and K.W. conducted the analyses, building the models, and prepared the figures. All authors reviewed the manuscript.
\subsection*{Additional information}
\vskip.4em
\noindent {\bf Supplementary information}: Available online along with all datasets used in the paper. It contains five real data applications, simulation results, further details on the methods and mathematical proofs.
\vskip.65em
\noindent {\bf \noindent Competing Interests}: The authors declare no competing interests.

\end{document}